\documentclass[10pt,a4paper,conference]{IEEEtran}

\usepackage[pdftex]{graphicx}
\usepackage[cmex10]{amsmath}
\usepackage{amsthm}
\usepackage{amssymb}
\usepackage{amscd}
\theoremstyle{plain}

\newtheorem{thm}{Theorem}
\newtheorem{prop}[thm]{Proposition}
\newtheorem{cor}[thm]{Corollary}
\newtheorem{lem}[thm]{Lemma}

\theoremstyle{definition}
\newtheorem{df}[thm]{Definition}

\newtheorem{eg}[thm]{Example}

\def\image{\textrm{image }}
\def\ed{\textrm{ed }}
\def\med{\textrm{med }}

\begin{document}
%
\title{The Nyquist theorem for cellular sheaves}

\author{\IEEEauthorblockN{Michael Robinson}
\IEEEauthorblockA{Department of Mathematics and Statistics\\
American University\\
4400 Massachusetts Ave NW\\
Washington, DC 20016\\
Email: michaelr@american.edu}}

\maketitle

\begin{abstract}
We develop a unified sampling theory based on sheaves and show that the Shannon-Nyquist theorem is a cohomological consequence of an exact sequence of sheaves.  Our theory indicates that there are additional cohomological obstructions for higher-dimensional sampling problems.  Using these obstructions, we also present conditions for perfect reconstruction of piecewise linear functions on graphs, a collection of non-bandlimited functions on topologically nontrivial domains.
\end{abstract}

\IEEEpeerreviewmaketitle

\section{Introduction}
\IEEEPARstart{T}{he} Shannon-Nyquist sampling theorem states that sampling a signal at twice its bandwidth is sufficient to reconstruct the signal.  Its wide applicability leads to the question of whether there exist similar conditions for reconstructing other data from samples in more general settings.  This article shows that perfect reconstruction for sampling of local algebraic data on simplicial complexes can be addressed through the machinery exact sequences of cellular sheaves.  As a demonstration of our technique, we recover the Nyquist theorem and generalize it to perfect reconstruction of piecewise linear signals on graphs.  Piecewise linear functions are not bandlimited, since their derivatives are not continuous.

\subsection{Historical context}

Sampling theory has a long and storied history, about which a number of recent survey articles \cite{benedetto_1990,feichtinger_1994,unser_2000,smale_2004} have been written.  Since sampling plays an important role in applications, substantial effort has been expended on practical algorithms.  Our approach is topologically-motivated, like the somewhat different approach of \cite{NiySmaWeiHom,chazal_2009}, so it is less constrained by specific timing constraints.  Relaxed timing constraints are an important feature of bandpass \cite{vaughan_1991} and multirate \cite{unser_1998} algorithms.  We focus on signals with local control, of which splines \cite{unser_1999} are an excellent example.

Sheaf theory has not been used in applications until fairly recently.  The catalyst for new applications was the technical tool of \emph{cellular sheaves}, developed in \cite{Shepard_1985}.  Since that time, an applied sheaf theory literature has emerged, for instance \cite{ghrist_2011,Lilius_1993,GhristCurryRobinson,RobinsonQGTopo,RobinsonLogic}.

Our sheaf-theoretic approach allows sufficient generality to treat sampling on non-Euclidean spaces.  Others have studied sampling on non-Euclidean spaces, for instance general Hilbert spaces \cite{pesenson_2001}, Riemann surfaces \cite{schuster_2004}, symmetric spaces \cite{ebata_2006}, the hyperbolic plane \cite{feichtinger_2011}, combinatorial graphs \cite{pesenson_2010}, and quantum graphs \cite{pesenson_2005,pesenson_2006}.  We show that sheaves provide unified sufficiency conditions for perfect reconstruction on abstract simplicial complexes, which encompass all of the above cases.

A large class of local signals are those with \emph{finite rate of innovation} \cite{groechening_1992, vetterli_2002}.  Our ambiguity sheaf is a generalization of the Strang-Fix conditions as identified in \cite{dragotti_2007}.  With our approach, one can additionally consider reconstruction using richer samples than simply convolutions with a function.

\section{Cellular sheaves}
\subsection{What is a sheaf?}
A sheaf is a mathematical object that stores locally-defined data over a space.  In order to formalize this concept, we need a concept of space that is convenient for computations.  The most efficient such definition is that of a simplicial complex.  

\label{sec:defs}
\begin{df}
An \emph{abstract simplicial complex} $X$ on a set $A$ is a collection of ordered subsets of $A$ that is closed under the operation of taking subsets.  We call each element of $X$ a \emph{face}.  A face with $k+1$ elements is called a $k$-face, though we usually call a $0$-face a \emph{vertex} and a $1$-face an \emph{edge}.  The \emph{face category} has as objects the elements of $X$ and as morphisms inclusions of one element of $X$ into another.
\end{df}

Although sheaves have been extensively studied over topological spaces (see \cite{Bredon} or the appendix of \cite{Hubbard} for a modern, standard treatment), the resulting definition is ill-suited for application to sampling.  Instead, we follow a substantially more combinatorial approach introduced in the 1985 thesis of Shepard \cite{Shepard_1985}.  

\begin{df}
A \emph{sheaf} $F$ on an abstract simplicial complex $X$ is a covariant functor from the face category of $X$ to the category of vector spaces.  Explicitly, 
\begin{itemize}
\item for each element $a$ of $X$, $F(a)$ is a vector space, called the \emph{stalk at $a$},
\item for each inclusion of two faces $a\to b$ of $X$, $F(a\to b)$ is a linear function from $F(a)\to F(b)$ called a \emph{restriction}, and 
\item for every composition of inclusions $a\to b \to c$, $F(b\to c) \circ F(a\to b)=F(a\to b\to c)$.
\end{itemize}
\end{df}

\begin{df}
Suppose $F$ is a sheaf on an abstract simplicial complex $X$ and that $\mathcal{U}$ is a collection of faces of $X$.  An assignment $s$ which assigns an element of $F(u)$ to each face $u\in\mathcal{U}$ is called a \emph{section} supported on $\mathcal{U}$ when for each inclusion $a\to b$ (in $X$) of objects in $\mathcal{U}$, $F(a\to b)s(a)=s(b)$.  A \emph{global section} is a section supported on $X$.  If $r$ and $s$ are sections supported on $\mathcal{U} \subset \mathcal{V}$, respectively, in which $r(a)=s(a)$ for each $a \in \mathcal{U}$ we say that \emph{$s$ extends $r$}.  The collection of sections supported on a given set forms a vector space.  
\end{df}

\begin{eg}
Consider $Y\subseteq X$ a subset of the vertices of an abstract simplicial complex. The functor $S$ which assigns a vector space $V$ to vertices in $Y$ and the trivial vector space to every other face is called a \emph{$V$-sampling sheaf supported on $Y$}.  To every inclusion between faces of different dimension, $S$ will assign the zero function.   For a finite abstract simplicial complex $X$, the space of global sections of a $V$-sampling sheaf supported on $Y$ is isomorphic to $\bigoplus_{y\in Y} V$.
\end{eg}

Recall that an abstract simplicial complex $X$ consists of \emph{ordered} sets.  For $a$ a $k$-face and $b$ a $k+1$-face, define 
\begin{equation*}
[b:a]=
\begin{cases}
+1&\text{if the order of elements in }a\text{ and }b\text{ agrees,}\\
-1&\text{if it disagrees, or}\\
0&\text{if }a\text{ is not a face of }b.\\
\end{cases}
\end{equation*}

\begin{eg}
Suppose $G$ is a graph in which each vertex has finite degree (evidently $G$ can be realized as an abstract simplicial complex).  Let $PL$ be the sheaf constructed on $G$ that assigns $PL(v)=\mathbb{R}^{1+\deg v}$ to each edge $v$ of degree $\deg v$ and $PL(e)=\mathbb{R}^2$ to each edge $e$.  The stalks of $PL$ specify the value of the function (denoted $y$ below) at each face and the slopes of the function on the edges (denoted $m_1,...,m_k$ below).  To each inclusion of a degree $k$ vertex $v$ into an edge $e$, let $PL$ assign the linear function $(y,m_1,...,m_e,...,m_k)\mapsto(y+[e:v]\frac{1}{2}m_e,m_e)$.  The global sections of this sheaf are \emph{piecewise linear functions} on $G$.  
\end{eg}

\begin{df}
A \emph{sheaf morphism} is a natural transformation between sheaves.  Explicitly, a morphism $f:F\to G$ of sheaves on an abstract simplicial complex $X$ assigns a linear map $f_a:F(a) \to G(a)$ to each face $a$ so that for every inclusion $a\to b$ in the face category of $X$, $f_b \circ F(a\to b) = G(a\to b) \circ f_a$.
\end{df}

\subsection{Sheaf cohomology}
Much of the theory of sheaves is concerned with computing spaces of sections and identifying obstructions to extending sections.  The machinery of cohomology systematizes the computation of the space of global sections for a sheaf.  

Define the following formal \emph{cochain} vector spaces $C^k(X;F)=\bigoplus_{a\text{ a }k\text{-face of }X}F(a)$.  The \emph{coboundary map} $d^k:C^k(X;F)\to C^{k+1}(X;F)$ takes an assignment $s$ from the $k$ faces to an assignment $d^ks$ whose value at a $k+1$ face $b$ is
\begin{equation*}
(d^ks)(b)=\sum_{a\text{ a }k\text{-face of} X} [b:a]F(a\to b)s(a).
\end{equation*}
It can be shown that $d^k\circ d^{k-1}=0$, so that the image of $d^{k-1}$ is a subspace of the kernel of $d^k$.

\begin{df}
The $k$-th \emph{sheaf cohomology} of $F$ on an abstract simplicial complex $X$ is 
\begin{equation*}
H^k(X;F)=\ker d^k / \image d^{k-1}.
\end{equation*}
\end{df}

Observe that $H^0(X;F)=\ker d^0$ consists precisely of those assignments $s$ which are global sections.  Cohomology is also a functor: sheaf morphisms induce linear functions between cohomologies.  This indicates that cohomology preserves and reflects the underlying relationships between sheaves.

\section{The Nyquist criterion for sheaves}
Suppose that $F$ is a sheaf on an abstract simplicial complex $X$, and that $S$ is a $V$-sampling sheaf on $X$ supported on a closed subcomplex $Y$.  A \emph{sampling} of $F$ is a morphism $s:F\to S$ that is surjective on every stalk.  Given a sampling, we can construct the \emph{ambiguity sheaf} $A$ in which the stalk $A(a)$ for a face $a\in X$ is given by the kernel of the map $F(a)\to S(a)$.  If $a\to b$ is an inclusion of faces in $X$, then $A(a\to b)$ is $F(a\to b)$ restricted to $A(a)$.  This implies that 
\begin{equation*}
\begin{CD}
0\to A \hookrightarrow F @>s>> S \to 0
\end{CD}
\end{equation*}
is an exact sequence, which induces the long exact sequence (via the Snake lemma)
\begin{equation*}
0\to H^0(X;A) \to H^0(X;F) \to H^0(X;S) \to H^1(X;A) \to
\end{equation*}
An immediate consequence is therefore
\begin{cor} (Sheaf-theoretic Nyquist theorem)
\label{thm:nyquist}
The global sections of $F$ are identical with the global sections of $S$\\ if and only if $H^k(X;A)=0$ for $k=0$ and $1$.
\end{cor}

The cohomology space $H^0(X;A)$ characterizes the \emph{ambiguity} in the sampling, while $H^1(X;A)$ characterizes its \emph{redundancy}.  Optimal sampling therefore consists of identifying minimal closed subcomplexes $Y$ so the resulting ambiguity sheaf $A$ has $H^0(X;A)=H^1(X;A)=0$.

Let us place bounds on the cohomologies of the ambiguity sheaf.  For a closed subcomplex $Y$ of $X$, let $F^Y$ be the sheaf whose stalks are the stalks of $F$ on $Y$ and zero elsewhere, and whose restrictions are either those of $F$ on $Y$ or zero as appropriate.  There is a surjective sheaf morphism $F\to F^Y$ and an induced ambiguity sheaf $F_Y$ which can be constructed in exactly the same way as $A$ before.  Thus, the dimension of each stalk of $F^Y$ is at least as large as that of any sampling sheaf, and the dimension of stalks of $F_Y$ are therefore as small as or smaller than that of any ambiguity sheaf.

\begin{prop} (Oversampling theorem)
\label{prop:oversample}
If $X^k$ is the closed subcomplex generated by the $k$-faces of $X$, then $H^k(X^{k+1};F_{X^k})=0$.
\end{prop} 
 
\begin{IEEEproof}
By direct computation, the $k$-cochains of $F_{X^k}$ are 
\begin{eqnarray*}
C^k(X^{k+1};F_{X^k})&=&C^k(X^{k+1};F)/C^k(X^k;F)\\
&=&\bigoplus_{a\text{ a }k\text{-face of}X} F(a) / \bigoplus_{a\text{ a }k\text{-face of}X} F(a)\\
&=& 0.
\end{eqnarray*}
\end{IEEEproof}

As an immediate consequence, $H^0(X;F_Y)=0$ when $Y$ is the set of vertices of $X$.

\begin{thm} (Sampling obstruction theorem)
\label{thm:obstruction}
Suppose that $Y$ is a closed subcomplex of $X$ and $s:F\to S$ is a sampling of sheaves on $X$ supported on $Y$.  If $H^0(X,F_Y) \not= 0$, then the induced map $H^0(X;F)\to H^0(X;S)$ is not injective.
\end{thm}
Succinctly, $H^0(X,F_Y)$ is an obstruction to the recovery of global sections of $F$ from its samples.

\begin{IEEEproof}
We begin by constructing the ambiguity sheaf $A$ as before so that
\begin{equation*}
\begin{CD}
0\to A \to F @>s>> S \to 0
\end{CD}
\end{equation*}
is a short exact sequence.  Observe that $S\to F^Y$ can be chosen to be injective, because the stalks of $S$ have dimension not more than the dimension of $F$ (and hence $F^Y$ also).  Thus the induced map $H^0(X;S)\to H^0(X;F^Y)$ is also injective.  Therefore, by a diagram chase on
\begin{equation*}
\begin{CD}
0\to H^0(X;A) @>>> H^0(X;F) @>s>> H^0(X;S)\\
&& @VV\cong V @VVV\\
0\to H^0(X;F_Y) @>>> H^0(X;F) @>>> H^0(X;F^Y)\\
\end{CD}
\end{equation*}
we infer that there is a surjection $H^0(X;A)\to H^0(X;F_Y)$.  By hypothesis, this means that $H^0(X;A)\not= 0$, so in particular $H^0(X;F)\to H^0(X;S)$ cannot be injective.
\end{IEEEproof}

\section{Applications}
\subsection{Bandlimited signals on the real line}
\label{sec:bandlimited}
\begin{figure}
\begin{center}
\includegraphics[width=3.25in]{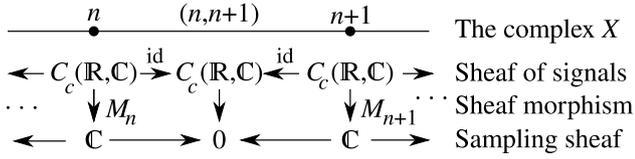}
\caption{The sheaves used in proving the traditional Nyquist theorem}
\label{fig:nyquist_setup}
\end{center}
\end{figure}

In this section, we prove the traditional form of the Nyquist theorem by showing that bandlimiting is a sufficient condition for $H^0(X;A)=0$.  We begin by specifying the following 1-dimensional simplicial complex $X$.  Let $X^0=\mathbb{Z}$ and $X^1=\{(n,n+1)\}$.  We construct the sheaf $C$ of signals (see Figure \ref{fig:nyquist_setup}) so that for every simplex, the stalk of $C$ is $C_c(\mathbb{R},\mathbb{C})$, the set of compactly supported complex-valued continuous functions, and each restriction is the identity.  Observe that the space of global sections of $C$ is therefore just $C_c(\mathbb{R},\mathbb{C})$.

Construct the sampling sheaf $S$ whose stalk on each vertex is $\mathbb{C}$ and each edge stalk is zero.  We construct a sampling morphism by the zero map on each edge, and by the inverse Fourier transform below on vertex $\{n\}$
\begin{equation*}
M_n(f)=\int_{-\infty}^\infty f(\omega) e^{-2\pi i n \omega} d\omega.
\end{equation*}
Then the ambiguity sheaf $A$ has stalks $C_c(\mathbb{R},\mathbb{C})$ on each edge, and $\{f\in C_c(\mathbb{R},\mathbb{C}) : M_n(f)=0\}$ on each vertex $\{n\}$.  

\begin{thm}(Traditional Nyquist theorem)
Suppose we replace $C_c(\mathbb{R},\mathbb{C})$ with the set of continuous functions supported on $[-B,B]$.  Then if $B\le 1/2$, the resulting ambiguity sheaf $A$ has $H^0(X;A)=0$.  Therefore, each such function can be recovered uniquely from its samples on $\mathbb{Z}$.
\end{thm}
\begin{IEEEproof}
The elements of $H^0(X;A)$ are given by the compactly supported continuous functions $f$ on $[-B,B]$ for which
\begin{equation*}
\int_{-B}^B f(\omega) e^{-2\pi i n \omega} d\omega=0
\end{equation*}
for all $n$.  Observe that if $B\le 1/2$, this is precisely the statement that the Fourier series coefficients of $f$ all vanish; hence $f$ must vanish.  This means that the only global section of $A$ is the zero function.  (Ambiguities can arise if $B>1/2$, because the set of functions $\{e^{-2\pi in\omega}\}_{n\in\mathbb{Z}}$ is then \emph{not} complete.)
\end{IEEEproof}

\subsection{Beyond Nyquist: Piecewise linear functions on graphs}
The sheaf-theoretic Nyquist theorem can treat nontrivial base space topologies as well as samples of different dimensions.  Consider the example of the sheaf of piecewise linear functions $PL$ on a graph, introduced in Section \ref{sec:defs} and the sampling morphism $s:PL\to PL^Y$ where $Y$ is a subset of the vertices of $X$.  Excluding one or two vertices from $Y$ does not prevent reconstruction in this case, because the samples include information about slopes along adjacent edges.

\begin{figure}
\begin{center}
\includegraphics[width=1.75in]{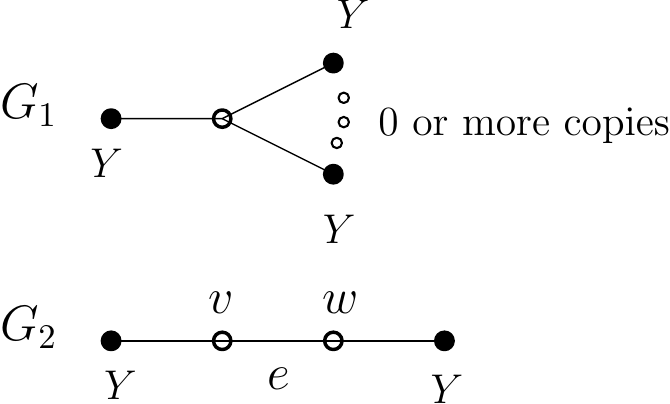}
\includegraphics[width=1.5in]{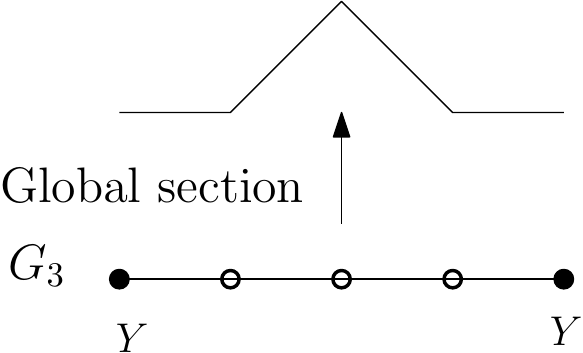}
\caption{Graphs $G_1$, and $G_2$ (left) and $G_3$ (right) for Lemma \ref{lem:pl_cases}.  Filled vertices represent elements of $Y$, empty ones are in the complement of $Y$.}
\label{fig:pl_cases}
\end{center}
\end{figure}

\begin{lem}
\label{lem:pl_cases}
Consider $PL_Y$, the subsheaf of $PL$ whose sections vanish on a vertex set $Y$ and the graphs $G_1$, $G_2$, and $G_3$ as shown in Figure \ref{fig:pl_cases}.  There are no nontrivial sections of $PL_Y$ on $G_1$ and $G_2$, but there are nontrivial sections of $PL_Y$ on $G_3$.
\end{lem}
\begin{IEEEproof}
If a section of $PL$ vanishes at a vertex $x$ with degree $n$, this means that the value of the section there is an $(n+1)$-dimensional zero vector.  The value of the section on every edge adjacent to $x$ is then the $2$-dimensional zero vector.  Since the dimensions in each stalk of $PL$ represent the value of the piecewise linear function and its slopes, linear extrapolation to the center vertex in $G_1$ implies that its value is zero too.

A similar idea applies in the case of $G_2$.  The stalk at $v$ has dimension 3.  Any section at $v$ that extends to the left must actually lie in the subspace spanned by $(0,0,1)$ (coordinates represent the value, left slope, right slope respectively).  In the same way, any section at $w$ that extends to the right must lie in the subspace spanned by $(0,1,0)$.  Any global section must extend to $e$, which must therefore have zero slope and zero value.

Finally $G_3$ has nontrivial global sections, spanned by the one shown in Figure \ref{fig:pl_cases}.
\end{IEEEproof}

\begin{df}
On a graph $G$, define the \emph{edge distance} between two vertices $v,w$ to be 
\begin{equation*}
\ed(v,w)=\begin{cases}
\min_p\{\text{\# edges in }p\text{ such that }p\text{ is a }\\\text{PL-continuous path from }v\to w\}\\
\infty\text{ if no such path exists}
\end{cases}
\end{equation*}
From this, the maximal distance to a vertex set $Y$ is 
\begin{equation*}
\med(Y)=\max_{x\in X^0} \{\min_{y\in Y}\; \ed(x,y)\}.
\end{equation*}
\end{df}

\begin{prop}(Unambiguous sampling)
Consider the sheaf $PL$ on a graph $X$ and $Y\subseteq X^0$.  Then $H^0(X;F_Y)=0$ if and only if $\med(Y) \le 1$.
\end{prop}

\begin{figure}
\begin{center}
\includegraphics[width=2.5in]{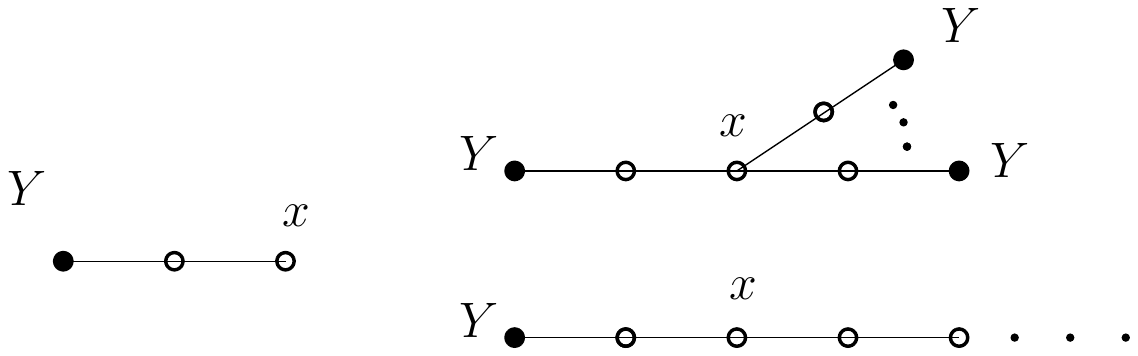}
\caption{The three families of subgraphs that arise when $\med(Y)>1$.  Filled vertices represent elements of $Y$, empty ones are in the complement of $Y$.}
\label{fig:pl_sampling}
\end{center}
\end{figure}

\begin{IEEEproof}
($\Leftarrow$) Suppose that $x\in X^0 \backslash Y$ is a vertex not in $Y$.  Then there exists a path with one edge connecting it to $Y$.  Whence we are in the case of $G_1$ of Lemma \ref{lem:pl_cases}, so any section at $x$ must vanish.

($\Rightarrow$) By contradiction.  Assume $\med(Y)>1$.  Without loss of generality, consider $x\in X^0\backslash Y$, whose distance to $Y$ is exactly 2.  Then one of the subgraphs shown in Figure \ref{fig:pl_sampling} must be present in $X$.  But case $G_3$ of Lemma \ref{lem:pl_cases} makes it clear that the most constrained of these (the middle panel of Figure \ref{fig:pl_sampling}) has nontrivial sections at $x$, merely looking at sections over the subgraph.
\end{IEEEproof}

\begin{prop}(Non-redundant sampling)
Consider the case of $s:PL\to PL^Y$.  If $Y=X^0$, then $H^1(X;A)\not= 0$.  If $Y$ is such that $\med(Y)\le 1$ and $|X^0\backslash Y|+\sum_{y\notin Y} \deg y = 2|X^1|$, then $H^1(X;A)=0$.
\end{prop}

\begin{IEEEproof}
The stalk of $A$ over each edge is $\mathbb{R}^2$, and the stalk over a vertex in $Y$ is trival.  However, the stalk over a vertex of degree $n$ not in $Y$ is $\mathbb{R}^{n+1}$.  Observe that if  $H^0(X;A)=0$, then $H^1(X;A)=C^1(X;A)/C^0(X;A)$.  Using the degree sum formula in graph theory, we compute that $H^1(X;A)$ has dimension $2|X^1| - \sum_{y\notin Y} (\deg y + 1)$.
\end{IEEEproof}

%
%

\section*{Acknowledgment}
This work was partly supported under Federal Contract No. FA9550-09-1-0643.


\bibliographystyle{IEEEtran}
\bibliography{sampta2013_bib}

\end{document}